\documentstyle[12pt]{amsart}

\def\bR{{\Bbb R}}
\def\bZ{{\Bbb Z}}

\def\cL{{\cal L}}
\def\cC{{\cal C}}
\def\cD{{\cal D}}
\def\Hom{{\cal{H}}om}

\def\om{\omega}

\def\pd#1{\frac{\partial}{\partial #1}}

\def\Vert{{\operatorname{Vert}}}
\newtheorem{thm}{Theorem}[section]
\newtheorem{lemma}[thm]{Lemma}
\newtheorem{prop}[thm]{Proposition}
\newtheorem{cor}[thm]{Corollary}
\newtheorem{defi}[thm]{Definition}
\newtheorem{rk}[thm]{Remark}

\def\prf{\noindent{\textsc{Proof}}\rm\ }
\def\endprf{\ \hfill $\Box$}
\begin{document}
\baselineskip=16pt
\parindent=0pt

\title[]{Hyper-symplectic structures on integrable systems}

\author[]{Claudio Bartocci and Igor Mencattini}
{ {\renewcommand{\thefootnote}{\fnsymbol{footnote}}
\footnotetext{\kern-15.3pt AMS Subject Classification: 53C26, 53D99, 37J35, 70H06.}
}}
 
\address{Universit\`a degli Studi di
Genova, Dipartimento di Matematica, Via Dodecaneso 35, 16146 Genova, ITALY}

\email{bartocci@@dima.unige.it}


\address{Boston University, Department of Mathematics and Statistics,
Boston University, 111 Cummington Street,
Boston, MA 02215, USA}

\email{igorre@@math.bu.edu}


\date{}
\maketitle

\begin{abstract} We prove that an integrable system over a
symplectic manifold whose symplectic form is covariantly
constant carries a natural hyper-symplectic structure.
Moreover, a special K\"ahler structure is induced on the base
manifold. 
\end{abstract}

\vspace{1cm}

\section{Introduction}
There exist deep relations between supersymmetric gauge theories and 
integrable systems, which were extensively investigated
by Donagi and Witten in their seminal paper \cite{DW}. In particular, 
it turns out that a key notion in both frameworks is
that of special K\"ahler manifold. According to the terminolgy fixed 
in \cite{F}, a special K\"ahler manifold is a K\"ahler
(or more generally, pseudo-K\"ahler) manifold equipped with a flat 
torsion-free connection $\nabla$ such that  the
covariant derivative of the K\"ahler form and the covariant 
differential of the complex structure are both equal to zero (for
details see Section 3). The cotangent bundle of any special K\"ahler 
manifold carries a compatible
hyper-symplectic structure; morevover, one can induce on it, by means 
of the connection $\nabla$, a natural structure of
algebraically completely integrable system \cite{F}.

In this note, we want to take an alternative point of view, looking 
at things the other way round. Indeed, our starting point
is a classical integrable system {\it \`a la}\ Liouville-Arnol'd, 
i.e.~a fibration $X\to B$, where $X$ is a symplectic manifold
and the fibres are Lagrangian tori. On such a fibration a flat 
torsion-free connection $\nabla$ is naturally defined. Assuming
the existence of a covariantly constant symplectic structure $\Omega$ 
on $B$, we prove the following:
\begin{enumerate}
\item $X$ admits a compatible hyper-symplectic structure (Theorem \ref{main});
\item $B$ is a special K\"ahler manifold (Corollary \ref{pippa}).
\end{enumerate}

Our approach should be compared to the treatment of special K\"ahler 
geometry developped by Hitchin
to describe the geometry of the moduli space $M$ of deformations of a 
complex Lagrangian submanifold of a complex K\"ahlerian
manifold \cite{H}. In particular, our construction implies the existence of a bi-Lagrangian immersion of the base manifold
$B$ into the total space $X$ (see Remark \ref{pippu}).

\section{The basic construction}

Let $(X,\om)$ be a connected symplectic manifold of dimension $2N$,
together with  a projection $\pi : X \to B$ whose
fibres $F_b=\pi^{-1}(b)$ are compact connected Lagrangian
submanifolds of $X$. We shall call $\pi : X \to B$ an integrable
system. By the Liouville-Arnol'd theorem
\cite{A}, the fibres
$F_b$ are isomorphic to $N$-dimensional tori. The homology groups
$H_1(F_b,\bZ)$ define a sheaf $\cL$
on the base manifold $B$, called the period lattice of the fibration
$\pi : X \to B$ \cite{D}. The dual sheaf of $\cL$,
$\Hom(\cL,\bZ)$ is isomorphic to the sheaf $R^1\pi_\ast\bZ$. The
sheaf of $\cC_B$- modules
$$\cD := \cL \otimes_{\bZ} \cC^{\infty}_B $$
is locally free; $\cL_b$ is a lattice in the fibre $D_b$ of the
associated vector bundle $D$ over $B$, and the quotient
$D_b/\cL_b$ is precisely the Liouville-Arnol'd torus at the point
$b$. We have the natural isomorphism $\pi^{\ast} D \simeq
\Vert (TX)$. We can define (locally) action-angle canonical
coordinates $(I_1,\dots,I_N,\varphi_1,\dots,\varphi_N)$: the
vertical Hamiltonian vector fields $X_i$ associated to the actions
$I_i$ single out a  basis $\{\gamma_i\}$ of local sections
of $\cL$ and the angle variables
$\{\varphi_i\}$ on the fibre $F_b$ satisfy the relation
$$\frac{1}{2\pi} \int_{\gamma_i} d\varphi_k = \delta^i_k\,.$$
The coordinates $I_i$ can be thought as affine coordinates on the
base manifold $B$; by associating $dI_i$ to
$X_i$, we obtain the natural isomorphism
\begin{equation}\label{ident}
T^\ast B \simeq D
\end{equation}
or equivalently $TB \simeq R^1\pi_\ast\bR \otimes_{\bR} \cC^{\infty}_B$.
The Gau\ss-Manin connection of the local system $R^1\pi_\ast\bR$
induces a torsion-free flat connection $\nabla$ on $TB$;
clearly, the vector fields $\partial/\partial I_i$ are parallel
w.r.t. $\nabla$. The subsheaf $\cL \subset \cD$ defines, via
the identification (\ref{ident}), a Lagrangian covering of $B$, that will be denoted by
$\Lambda \subset T^\ast B$; the  monodromy of $\Lambda$ coincides with
the holonomy of $\nabla$. According to \cite{D}, if the fibration
$\pi : X \to B$ admits a section, then there is an
isomorphism $X \simeq T^\ast B/ \Lambda$ (which is an isomorphism of
symplectic manifolds fibred in Lagrangian tori if the
section is Lagrangian).

One has the identification $\Vert (TX) \simeq \pi^\ast T^\ast B$. Let
us consider the dual fibration
$\hat X = R^1\pi_\ast\bR/ R^1\pi_\ast\bZ$, naturally endowed with a
projection $\hat\pi:
\hat X \to B$.  The identification $\Vert (T\hat X) \simeq \hat\pi^\ast T
B$ can be plugged in the Atiyah sequence
\begin{equation}
0 \to \Vert (T\hat X) \to T\hat X \to \hat\pi^\ast T B \to 0
\end{equation}
which is splitted by the Gau\ss-Manin connection. In this way, we get
a decomposition
\begin{equation}
T\hat X = \hat\pi^\ast T^\ast B \oplus \hat\pi^\ast T B\,.
\end{equation}
According to \cite{BMP,AP}, an almost-complex structure $J$ is defined  on
$\hat X$ by setting $J(U,V)=(-V,U)$; actually, $J$ is integrable
because the Gau\ss-Manin connection is flat. If we introduce on $\hat
X$ local coordinates $(I_1,\dots,I_N,
\Phi_1,\dots,\Phi_N)$, where the $\Phi$'s are dual coordinates to the
$\varphi$'s, then the holomorphic coordinates are given
by $z_j = I_j +
i \Phi_j$. If the fibration $\pi : X \to B$ has a section, then $X$
can be identified to $\hat X$, in such a way that the
natural lifting of the Gau\ss-Manin connection coincides with the
connection we have introduced on $\hat X$ and the $\Phi$'s
coincide with the angles $\varphi$'s.

Let us suppose that the base manifold $B$ has a sympletic structure
$\Omega$; this implies $N=2n$. We shall be concerned only
with local properties, so that we can assume that a Lagrangian
section $\sigma_{\omega}: B \to X$ does exist.

If $\Omega$ is covariantly costant, i.e. $\nabla \Omega=0$, then, the
splitting $T X = \pi^\ast T^\ast
   B \oplus \pi^\ast T^\ast B$ (together with the isomorphsim $T
B\simeq T^\ast B$ provided by $\Omega$ itself) can be used to
construct a new symplectic form $\chi= -\Omega\oplus \Omega$ on $X$.
The condition $\nabla\Omega=0$ is equivalent to the fact
that the action coordinates $I_i$ can be assumed to be canonical
w.r.t. $\Omega$; dropping the indices, we shall write
short
$I= (x,y)$ and
$\Omega = dx\wedge dy$. By using the same notations, we also have
$\varphi=(p,q)$ and $\omega= d\varphi \wedge dI = dp\wedge
dx + dq\wedge dy$. Finally, we get
$$\chi = - dp \wedge dq + dx \wedge dy\,.$$

We shall denote by $J_\omega$ and $J_\chi$ the complex structures
associated, respectively, to the symplectic structures
$\omega$ and $\chi$
   on the dual fibration $\hat X$ (in principle, the complex structure
$J_\chi$ is defined only locally). In the
local coordinates $(x,y,P,Q)$ on $\hat X$, one has:
\begin{align*}
J_\omega &= dx \otimes  \pd{P} - dP\otimes \pd{x} + dy\otimes \pd{Q}
- dQ\otimes \pd{y} \\
J_\chi &= - dP\otimes \pd{Q} + dQ \otimes \pd{P} + dx \otimes \pd{y}
- dy \otimes \pd{x}\,.
\end{align*}
Let us define $K= J_\omega\circ J_\chi$.
\begin{lemma}\label{lemma} $K$ is a complex structure.
\end{lemma}
\prf An easy computation yields:
$$\begin{array}{r@{\,;\quad}l} K (dx) = dQ & K (dy) = -dP \\
   K(dQ) = -dx & K(dP) = dy\,.
\end{array}$$
Therefore, $K\circ K = -\operatorname{Id}$.
\endprf

In local coordinates, we can write: 
$$K = - dx \otimes \pd{Q} +dQ\otimes\pd{x} - dP\otimes\pd{y} + dy\otimes \pd{P}.$$
The holomorphic coordinates
$(\alpha_1,\dots\alpha_n,\beta_1,\dots,\beta_n)$ induced by $K$
satisfy the relations $d\alpha= dx + i
dQ$, $d\beta = dP + i dy$.  On $X$ we define the symplectic form
$$\sigma = dq \wedge dx + dy \wedge dp\,.$$
We shall write $J_\sigma$ instead of  $K$. Summing up, we have proved
the following result.
\begin{thm}\label{main} Let $X\to B$ be an integrable system over a
symplectic manifold $(B,\Omega)$ such that $\nabla
\Omega=0$. Then there exist (locally) a hyper-symplectic structure on
$X$  and a hyper-complex structure on the
dual fibration $\hat X$.
\end{thm}
\prf Lemma \ref{lemma} shows that $J_\omega$, $J_\chi$ and $J_\sigma$
determine a hypercomplex structure on $\hat X$. It
follows that $(\omega^{-1} \chi)^2 = -\operatorname{Id}_{TX}$; hence,
$\omega$, $\chi$ and $\sigma$ define a
hyper-symplectic structure on
$X$ \cite{X}. \endprf

We can identify $X$ and $\hat X$ in such a way that the coordinates
$(x,y,p,q)$ correspond to the coordinates $(x,y,P,Q)$.
As a straighforward consequence of Theorem \ref{main} we obtain
the following
\begin{cor}\label{corollario}
Let $\rho_\omega: B\to X$ be a Lagrangian immersion w.r.t.~the
symplectic structure $\omega$. Then the image $\rho_\omega(B)$
is a complex submanifold w.r.t.~the complex structure $J_\chi$. \endprf
\end{cor}
Since the holomorphic coordinates on $X$ induced by $J_\chi$ are
$(u,v)=(q+ip, x+iy)$, we can consider also $q+ip$ as
homolomorphic coordinates on $B$ (to be precise, we need to identify
$TB$ and $\Vert(TX)$ by means of the symplectic form
$\omega$). Given a Lagrangian immersion $\rho_\sigma$ (resp.,
$\rho_\chi$) w.r.t. $\sigma$ (resp., $\chi$), the analogue of
Corollary
\ref{corollario} is readily stated:  $\rho_\sigma(B)$
is a complex submanifold w.r.t.~$J_\omega$ (resp., $\rho_\sigma(B)$
is a complex submanifold w.r.t.~$J_\omega$).

\section{Special K\"ahler geometry on $B$}

In this section we show how to recover from the data described in the 
previous sections a special K\"ahler structure on $B$
\cite{F, H}. We start with the following basic definition.
\begin{defi}\cite{Co} A complex manifold $(B,I)$ is called {\it 
special complex} if there is a flat, torsion-free
(linear) connection $\nabla$ such that:
$$d_{\nabla} I=0\,.$$
If there is a covariantly constant symplectic form $\Omega$ on $B$ 
(i.e.~$\nabla\Omega=0$), then the triple $(M,I,\Omega)$
is called {\it special symplectic.}
\end{defi}

Let us suppose that a Lagrangian section 
$\rho_{\sigma}:B\longrightarrow X$ does exist. From Corollary 
\ref{corollario}, we
have that $\rho_{\sigma}(B)$ is a complex submanifold w.r.t.~the 
complex structure $J_\omega$. Moreover the holomorphic
coordinates induced by this complex structure on $X$ are 
$(z,w)=(x+ip,y+iq)$. These coordinates define, via
$\rho_{\sigma}$, a double set of holomorphic coordinates on $B$, 
$z=x+ip$ and $w=y+iq$, while the restriction to
$\rho_{\sigma}(B)$ of the complex structure $J_{\omega}$ induces a 
complex structure on $B$ that can be written in local
coordinates as:
\begin{equation}
I=- \big(dp\otimes \pd{x}+ dq\otimes \pd{y}\big)\label{com}.
\end{equation}
Since the coordinates $(x,y)$ are flat w.r.t the Gau\ss-Manin 
connection we have:
$$d_{\nabla} I=0.$$
So have proven the following result.
\begin{prop}\label{pippo} Under the same assumptions of Theorem 
\ref{main} and assuming the existence of a Lagrangian
section as above, the manifold $B$ is special symplectic.\endprf
\end{prop}

Let us recall that a special symplectic manifold 
$(B,\Omega,I)$ is said to be {\it special K\"ahler} if the
symplectic form $\Omega$ is $I$-invariant \cite{Co, F, H}.
Notice that the 2-tensor  $g(\cdot,\cdot)= \Omega(\cdot,I\cdot)$, whilst is always symmetric, in general is not positive 
definite. So, the name of special pseudo-K\"ahler would be more appropriate.

The symplectic form $\Omega$ can be written as $\Omega=dx\wedge dy$ 
(recall that the $(x,y)$'s are the flat
symplectic coordinates on
$B$) and a  straighforward calculation shows that such form is 
invariant w.r.t.~the complex structure defined in (\ref{com}). Thus,
  under the same assumption as in Proposition \ref{pippo}, the 
following result is easily proved.
\begin{cor}\label{pippa} The manifold $B$ is special K\"ahler, with 
(pseudo)-K\"ahler metric $g$ given by $g=\Omega\circ I$.
\end{cor}

\begin{rk}\label{pippu}\rm
In \cite{H}, Theorem 2.4, Hitchin characterizes special (pseudo) 
K\"ahler manifolds as manifolds that can be (locally) identified as 
bi-Lagrangian
submanifolds of $V\times V$, where $V$ is a symplectic vector space. 
In our description of  the special geometry defined on the base of 
the fibration,
we use the existence of only one Lagrangian section $\rho_{\sigma}$ 
even if $B$ can be (clearly) identified via the  Gau\ss-Manin 
connection with a $\omega$-Lagrangian submanifold of the total space $X$.
\end{rk}

\vskip20pt
\paragraph{\bf Acknowledgements.}
C.B. acknowledges the financial support of the MIUR and the University
of Genova through the national research project ``Geometria dei 
sistemi integrabili''. The authors thank Ugo Bruzzo
for valuable discussions.

\end{document}